\input amstex
\magnification=\magstep1 
\baselineskip=13pt
\documentstyle{amsppt}
\vsize=8.7truein \CenteredTagsOnSplits \NoRunningHeads
\def\conv{\operatorname{conv}}
\def\spa{\operatorname{span}}
\def\rk{\operatorname{rank}}
\def\rkp{\operatorname{rank_{psd}}}
\def\sym{\operatorname{Sym}}
 
 \topmatter

\title  Approximations of convex bodies by polytopes and by projections of spectrahedra  \endtitle
\author Alexander Barvinok  \endauthor
\address Department of Mathematics, University of Michigan, Ann Arbor,
MI 48109-1043, USA \endaddress
\email barvinok$\@$umich.edu \endemail
\date April 2012 \enddate
\thanks  This research was partially supported by NSF Grant DMS 0856640.
\endthanks 
\abstract We prove that for any compact set $B \subset {\Bbb R}^d$ and for any $\epsilon >0$ there is 
a finite subset $X \subset B$ of $|X|=d^{O(1/\epsilon^2)}$ points such that the maximum absolute value of any linear 
function $\ell: {\Bbb R}^d \longrightarrow {\Bbb R}$ on $X$ approximates the maximum absolute value of $\ell$ on 
$B$ within a factor of $\epsilon \sqrt{d}$. We also discuss approximations of convex bodies by projections of spectrahedra,
that is, by projections of sections of the cone of positive semidefinite matrices by affine subspaces.
\endabstract
\endtopmatter

\document

\head 1. Introduction and main results \endhead

We present two results on approximating general convex bodies by efficiently computable convex bodies in the general spirit
of \cite{BV08}. Having fixed a compact set $B \subset {\Bbb R}^d$, we are interested in constructing an algorithm, preferably
of a reasonably low complexity, which allows us to approximate the maximum of a given linear function 
$\ell: {\Bbb R}^d \longrightarrow {\Bbb R}$ on $B$. 

Our first result describes how well we can approximate $B$ by a finite subset $X \subset B$ of a controlled size.

\proclaim{(1.1) Theorem} Let $B \subset {\Bbb R}^d$ be a compact set. Then for any positive integer $k$ there exists a set $X \subset B$ such that for the cardinality $|X|$ of $X$ we have
$$|X| \ \leq \ 1 + {1 \over 2} {d+k-1 \choose  k} +{1 \over 2}{d+k-1 \choose k}^2$$ and  such that
$$\max_{x \in B} |\ell(x)| \ \leq \  {d+k-1 \choose  k}^{1\over 2k} \max_{x \in X} |\ell(x)|$$
for every linear function $\ell: {\Bbb R}^d \longrightarrow {\Bbb R}$.
\endproclaim 

Let us fix $k$ and let the dimension $d$ grow. Then the cardinality $|X|=d^{O(k)}$ of $X$ is polynomial in $d$ while for the 
approximation factor we have
$${d +k -1 \choose k}^{1\over 2k} \approx {\sqrt{d} \over (k!)^{1/2k}} \approx \sqrt{d e \over k}.$$
In particular, choosing a sufficiently large $k$ we can replace $B$ by a set $X$ of polynomially many in $d$ points such that 
the maximum absolute value any linear function on $X$ approximates the maximum absolute value of the function on 
$B$ within a factor $\epsilon \sqrt{d}$ for any $\epsilon >0$, fixed in advance. To obtain a constant factor approximation 
we have to choose $X$ consisting of exponentially many in $d$ points.

Next, we consider approximations of convex bodies by more complicated sets.

In the space ${\Bbb R}^{r \times r}$ of $r \times r$ real matrices we consider the closed convex cone 
${\Bbb R}^{r \times r}_+$ of symmetric positive semidefinite matrices. A section of ${\Bbb R}^{r \times r}_+$ by an affine 
subspace is called sometimes a {\it spectrahedron}, see, for example, \cite{GN11}. The problem of optimizing a given 
linear function on the affine image (projection) of a spectrahedron is a problem of {\it semidefinite programming}, which, under some
technical qualifications, can 
be solved in polynomial time, see, for example, \cite{Tu10}. 

The following result was obtained by Gouveia, Parrilo and Thomas \cite{G+10}, \cite{GT10} in the language of 
{\it theta bodies}. Nevertheless, we give a proof here as it connects the topic with the concept of {\it positive semidefinite 
rank} of a matrix introduced in \cite{F+11} and \cite{G+11} and raises some interesting questions. We note, however, that the 
proof uses the same idea as the proof from \cite{GT10}, only stated in a different language.

\proclaim{(1.2) Theorem} Let $B \subset {\Bbb Z}^d$ be a finite set of integer vectors.
 For any positive integer
$k$ there exists a convex set $C \subset {\Bbb R}^d$ such that the following holds:
\roster 
\item We can write $C$ as a Minkowski sum $C=C' + {\Cal L}^{\bot}$, where ${\Cal L} \subset {\Bbb R}^d$ is a subspace 
and $C' \subset {\Cal L}$ is an affine image of a section of the cone of 
$r \times r$ symmetric positive semidefinite matrices by an affine subspace with 
$$r \ \leq \ {d +k+2 \choose k},$$ 
\item We have 
$$B \subset C$$
and
\item For any linear function 
$$\ell(x)=a_1 x_1 + \ldots +a_d x_d \quad \text{for} \quad x=\left(x_1, \ldots, x_d\right)$$
with integer coefficients $a_1, \ldots, a_d$ such that 
$$\max_{x \in B} \ell(x) -\min_{x \in B} \ell(x) \ \leq \ k$$
we have 
$$\max_{x \in B} \ell(x) = \max_{x \in C} \ell(x).$$
\endroster
\endproclaim

If $k$ is fixed in advance then $r=d^{O(k)}$ is bounded by a polynomial in the dimension and $C$ approximates the convex hull $\conv(B)$ precisely with respect to any lattice direction for which 
the width of the convex hull is bounded by $k$. If we allow the lattice width $k$ to be linear in $d$, the dimension $r$ of the 
ambient space for the spectrahedron becomes exponentially large in $d$.

The paper is organized as follows.

In Section 2, we discuss some preliminaries concerning (symmetric) tensor powers of spaces.

In Section 3, we prove Theorem 1.1.

In Section 4, we prove Theorem 1.2.

In Section 5, we present a related result on approximating non-negative matrices by matrices with a small 
{\it positive semidefinite rank}, studied in \cite{F+11} and \cite{G+11}.

\head 2. Preliminaries \endhead 

We consider Euclidean space ${\Bbb R}^d$ with scalar product 
$$\langle x, y \rangle = \sum_{i=1}^d x_i y_i \quad \text{for} \quad x=\left(x_1, \ldots, x_d\right) \quad \text{and} \quad 
y=\left(y_1, \ldots, y_d\right)$$
and the corresponding Euclidean norm
$$\|x\| =\sqrt{\langle x, x \rangle}.$$
For a positive integer $k$ we interpret the tensor product
$$\left({\Bbb R}^d\right)^{\otimes k} =\underbrace{{\Bbb R}^d \otimes \ldots \otimes {\Bbb R}^d}_{\text{$k$ times}}$$
as $d^k$-dimensional Euclidean space of arrays $X=\left(x_{i_1 \ldots i_k}\right)$, $1 \leq i_1, \ldots, i_k \leq d$, with 
scalar product 
$$\bigl\langle X, Y \big\rangle =\sum_{1 \leq i_1, \ldots, i_k \leq d} x_{i_1 \ldots i_k} y_{i_1 \ldots i_k} 
\quad \text{for} \quad X=\left(x_{i_1 \ldots i_k}\right) \quad \text{and} \quad Y=\left(y_{i_1 \ldots i_k} \right).$$
For a vector $x \in {\Bbb R}^d$, we define $x^{\otimes k} \in \left({\Bbb R}^d\right)^{\otimes k}$ by 
$$\left(x^{\otimes k}\right)_{i_1 \ldots i_k} =x_{i_1} \cdots x_{i_k} \quad \text{for} \quad 
1 \leq i_1, \ldots, i_k \leq d \quad \text{and} \quad x=\left(x_1, \ldots, x_d \right).$$
We have 
$$\big\langle x^{\otimes k},\ y^{\otimes k} \big\rangle =\langle x, y \rangle^k \quad \text{for all} \quad 
x, y \in {\Bbb R}^d.$$
We observe that $x^{\otimes k}$ lies in the subspace $\sym\left({\Bbb R}^d\right)^{\otimes k}$ consisting of the tensors 
$X=\left(x_{i_1\ldots i_k}\right)$ for which values of the coordinates do not change when indices $i_1, \ldots i_k$ are 
permuted. We have 
$$\dim \sym\left({\Bbb R}^d\right)^{\otimes k} ={d+k-1 \choose k}.$$
Finally, the space $\left({\Bbb R}^d\right)^{\otimes 2}$ is naturally identified with the space ${\Bbb R}^{d \times d}$ of $d \times d$ matrices while 
subspace $\sym \left({\Bbb R}^d\right)^{\otimes 2}$ is identified with the subspace of symmetric matrices. 
For $x \in {\Bbb R}^d$ the matrix $x^{\otimes 2}$ is positive semidefinite. As is well-known, the cone ${\Bbb R}^{d \times d}_+$ of symmetric positive semidefinite 
matrices is spanned by matrices $x^{\otimes 2}$. We have $\big\langle X, Y \big\rangle \geq 0$ for any $d \times d$ 
symmetric positive semidefinite matrices.

\head 3. Proof Theorem 1.1 \endhead

We start with a lemma.
\proclaim{(3.1) Lemma} Let $B \subset {\Bbb R}^d$ be a compact set. Then there exists a set $X \subset B$ 
of not more than $1+{d(d+1) \over 2}$ points such that for any linear function $\ell: {\Bbb R}^d \longrightarrow {\Bbb R}$
one has 
$$\max_{x \in B} |\ell(x)| \ \leq \ \sqrt{d} \max_{x \in X} |\ell(x)|.$$
\endproclaim
\demo{Proof} Without loss of generality we assume that $B$ spans ${\Bbb R}^d$. We consider the (necessarily unique) ellipsoid 
$E \subset {\Bbb R}^d$ centered at the origin and of the minimum volume among those which contain $B$, see, for example,
\cite{Ba97}.  Applying an invertible linear transformation,
if necessary, we assume that $E=\left\{x \in {\Bbb R}^d: \ \|x\| \leq 1 \right\}$ is the unit ball. F. John's conditions,
see \cite{Ba97}, 
state that there is a finite subset $X \subset B$ and numbers $\alpha_x \geq 0$ such that 
$$\aligned &\sum_{x \in X} \alpha_x x^{\otimes 2} ={1 \over d} I \\ &\qquad \text{and} \\
&\sum_{x \in X} \alpha_x =1 \endaligned \tag3.1.1$$
where $I$ is the $d \times d$ identity matrix. 

For completeness, we sketch a proof of  (3.1.1). 
If matrix $d^{-1} I$ does not lie in the convex hull of the compact set
$\left\{x^{\otimes 2}: \ x \in B\right\}$ then it can be separated from the set by an affine hyperplane, 
which implies that there is a quadratic form $q: {\Bbb R}^d \longrightarrow {\Bbb R}$ such that 
$q(x) \leq 1$ for all $x \in B$ and such that $\text{trace\ }q >d$. Then for a sufficiently small $\epsilon >0$
the ellipsoid $\tilde{E}$ defined by the inequality 
$$(1-\epsilon) \|x\|^2 +\epsilon q(x) \ \leq \ 1$$
contains $B$ and has a smaller volume, which is a contradiction.

Carath\'eodory's Theorem then implies that we can choose 
$|X| \ \leq 1+{d(d+1) \over 2}$. 
The first equation of (3.1.1) can be also written as 
$$\sum_{x \in X} \alpha_x \langle c, x \rangle^2 ={1 \over d} \|c\|^2 \quad \text{for all} \quad c \in {\Bbb R}^d,$$
from which it follows that 
$$\max_{x \in X} |\langle c, x \rangle| \ \geq \ {1 \over \sqrt{d}} \|c\| \ \geq \ {1 \over \sqrt{d}}\max_{x \in B} |\langle c, x \rangle|.$$
{\hfill \hfill \hfill} \qed
\enddemo

\subhead (3.2) Proof of Theorem 1.1 \endsubhead 
Let us consider
$$B_k =\left\{ x^{\otimes k}: \quad x \in B \right\}.$$
Thus $B_k$ is a compact subset of a ${d+k-1 \choose k}$-dimensional space $\sym \left({\Bbb R}^d\right)^{\otimes k}$. 
Applying Lemma 3.1 we conclude that there exists a subset $X \subset B$ such that 
$$|X| \ \leq \  1+{1 \over 2} {d+k-1 \choose k} +{1 \over 2} {d+k-1 \choose k}^2 $$ and such such that 
$$\max_{x \in B} \left| \big\langle Y,\ x^{\otimes k} \big\rangle \right| \ \leq \ 
{d+k-1 \choose k}^{1/2} \max_{x \in X} \left| \big\langle Y,\ x^{\otimes k} \big\rangle \right|$$
for any $Y \in \sym \left({\Bbb R}^d\right)^{\otimes k}$. Choosing $Y=y^{\otimes k}$ we conclude that 
$$\max_{x \in B} \left| \langle y, x \rangle \right| \ \leq \ 
{d+k-1 \choose k}^{1/2k} \max_{x \in X} \left| \langle y, x \rangle \right|$$
for any $y \in {\Bbb R}^d$, which completes the proof.
{\hfill \hfill \hfill} \qed

It follows from our proof that we can choose the set $X$ among the contact points of the ellipsoid of the minimum volume
centered at the origin and containing the set $B_k =\left\{ x^{\otimes k}: \quad x \in B \right\}$.

\head 4. Proof of Theorem 1.2 \endhead 

We deal with real matrices 
$A=\left(a_{ij}\right)$ for $i \in I$ and $j \in J$, where $I$ and $J$ are possibly infinite sets of indices. We say that 
$$\rk A \ \leq \ n$$
if there exist vectors $u_i \in {\Bbb R}^n$ for $i \in I$ and $v_j \in {\Bbb R}^n$ for $j \in J$ such that 
$$a_{ij} = \langle u_i, v_j \rangle \quad \text{for all} \quad i \in I \quad \text{and all} \quad  j \in J.$$
We say that 
$$\rk A=n$$ if $n$ is the smallest non-negative integer satisfying $\rk A \leq n$. We can define $\rk A=\infty$ if there 
is no such $n$, but we will only deal with matrices of a finite rank.
Our definition agrees with the usual definition of the rank of a matrix, when $I$ and $J$ are finite.

We need some concepts and results of \cite{F+11} and \cite{G+11}.

\definition{(4.1) Definition} Let $A=\left(a_{ij}\right)$, $i \in I$, $j \in J$ be a non-negative matrix. We say that 
$$\rkp A \ \leq \ r$$
if there exist $r \times r$ symmetric positive semidefinite matrices $U_i$ for $i \in I$ and $V_j$ for  
$j \in J$ such that 
$$a_{ij} =\big\langle U_i,\ V_j \big\rangle \quad \text{for all} \quad i, j.$$
\enddefinition 

The following result was proved in \cite{F+11} and \cite{G+11}. For completeness, we present its proof here.

\proclaim{(4.2) Lemma} Let $\left\{u_i:\  i \in I\right\} \subset {\Bbb R}^d$ and  
$ \left\{v_j:\  j \in J \right\} \subset {\Bbb R}^d$ be sets of vectors such that 
$$\langle u_i, v_j \rangle \ \leq \ 1 \quad \text{for all} \quad i, j.$$
Suppose further that $\spa\left(v_j:\ j \in J \right)={\Bbb R}^d$.

Let us define matrix $A=\left(a_{ij}\right)$ by 
$$a_{ij}=1 -\langle u_i, v_j \rangle \quad \text{for all} \quad i, j.$$
Suppose that 
$$\rkp A \ \leq \ r.$$
Then there exists a convex set $C \subset {\Bbb R}^d$ which is an affine image of a section of the cone of
$r \times r$ symmetric positive semidefinite matrices by an affine subspace such that
$$u_i \in C \quad \text{for all} \quad i \in I$$
and 
$$\langle x, v_j \rangle \ \leq \ 1 \quad \text{for all} \quad j \in J \quad \text{and all} \quad x \in C.$$
\endproclaim
\demo{Proof} Since $\rkp A \ \leq \ r$ there exist $r \times r$ positive semidefinite matrices $U_i$ and 
$V_j$ such that 
$$1-\langle u_i, v_j \rangle =\big\langle U_i, V_j \big\rangle \quad \text{for all} \quad i, j. \tag4.2.1$$
Let us define an affine subspace ${\Cal L} \subset {\Bbb R}^d \oplus {\Bbb R}^{r \times r}$ by the equations 
$$1- \langle x, v_j \rangle = \big\langle X, V_j \big\rangle \quad \text{for} \quad j \in J, \tag4.2.2$$
where $x \in {\Bbb R}^d$ and $X \in {\Bbb R}^{r \times r}$. The map 
$$(x, X) \longmapsto X$$
projects ${\Cal L}$ onto an affine subspace ${\Cal A} \subset {\Bbb R}^{r \times r}$. 
We claim that for every $X \in {\Cal A}$ there is a unique $x \in {\Bbb R}^d$ such that $(x, X) \in {\Cal L}$. Indeed, 
if $(x, X) \in {\Cal L}$ and $(y, X) \in {\Cal L}$, then 
$$\langle x-y,\ v_j \rangle =0 \quad \text{for} \quad j \in J$$
and since the set $\left\{v_j: j \in J\right\}$ spans ${\Bbb R}^d$ we conclude that $x=y$. This allows us to define an affine map (projection)
$$T: {\Cal A} \longrightarrow {\Bbb R}^d$$
by letting 
$$T(X)=x \quad \text{if} \quad (x, X) \in {\Cal L}.$$
We let 
$$C=T\left( {\Cal A} \cap {\Bbb R}^{r \times r}_+\right).$$
Since by (4.2.1)--(4.2.2) we have $\left(u_i,\ U_i\right) \in {\Cal L}$, we conclude that $U_i \in {\Cal A}$ and $T\left(U_i\right)=u_i$. Since 
$U_i \in {\Bbb R}^{r \times r}_+$, we have $u_i \in C$ for all $i \in I$. 

Let us pick any $x \in C$. Then there exists an $X \in {\Bbb R}^{r \times r}_+$ such that $(x, X)$ satisfies (4.4.2).
Since $V_j \in {\Bbb R}^{r \times r}_+$ we have $\langle X, V_j \rangle \geq 0$ for all $j$ and hence 
$\langle x, v_j \rangle \leq 1$ for all $j \in J$.
{\hfill \hfill \hfill} \qed
\enddemo

The following observation is also from \cite{F+11} and \cite{G+11}.

\proclaim{(4.3) Lemma} Let $A=\left(a_{ij}\right)$ and $B=\left(b_{ij}\right)$ for $i \in I$ and $j \in J$ be matrices such that 
$$a_{ij}=b_{ij}^2 \quad \text{for all} \quad i, j.$$
Then 
$$\rkp A \ \leq \ \rk B.$$
\endproclaim
\demo{Proof}
Let $\rk B=d$. Then there exist vectors $u_i, v_j \in {\Bbb R}^d$ such that 
$$b_{ij}=\langle u_i, v_j \rangle \quad \text{for all} \quad i \in I \quad \text{and all} \quad j \in J.$$
Then
$$a_{ij}=\big\langle u_i^{\otimes 2}, v_j^{\otimes 2} \big\rangle \quad \text{for all} \quad i, j.$$ 
Since $u_i^{\otimes 2}$ and $v_j^{\otimes 2}$ are $d \times d$ positive semidefinite matrices, the result follows.
{\hfill \hfill \hfill} \qed
\enddemo

The following result is a standard linear algebra fact.

\proclaim{(4.4) Lemma} Let $A=\left(a_{ij}\right)$ be a real matrix and let 
$p: {\Bbb R} \longrightarrow {\Bbb R}$ be a polynomial of degree $k$. Let us define a matrix $B=\left(b_{ij}\right)$ 
by 
$$b_{ij}=p\left(b_{ij}\right) \quad \text{for all} \quad i\in I \quad \text{and all} \quad j \in J.$$
Then 
$$\rk B \ \leq \ {k +\rk A \choose k}.$$
\endproclaim
\demo{Proof}
We write 
$$p(t)=\sum_{m=0}^k \alpha_m t^m$$
for some $\alpha_m \in {\Bbb R}$. 
Let $\rk A=d$, so
$$a_{ij} =\langle u_i, v_j \rangle \quad \text{for all} \quad i, j$$ 
and some vectors $u_i:  i \in I$ and  $v_j:  j \in J$ in ${\Bbb R}^d$.
Then we can write 
$$b_{ij}=p\left(a_{ij}\right) =\alpha_0 + \sum_{m=1}^{k} \alpha_j \langle u_i, v_j \rangle^m = 
\alpha_0 + \sum_{m=1}^{k} \alpha_m \big\langle   u_i^{\otimes m}, v_j^{\otimes m} \big\rangle.$$
Let us introduce vectors 
$$U_i =\alpha_0 \oplus \alpha_1 u_i \oplus \alpha_2 u_i^{\otimes 2}  \oplus \ldots \oplus \alpha_{k} u_i^{\otimes k} \quad \text{and} \quad 
V_j=1 \oplus v_j \oplus v_j^{\otimes 2} \oplus  \ldots \oplus v_j^{\otimes k}$$ in Euclidean space 
$${\Bbb R} \oplus {\Bbb R}^d \oplus \left({\Bbb R}^d\right)^{\otimes 2} \oplus \ldots \oplus \left({\Bbb R}^d\right)^{\otimes k}.$$
Hence we can write 
$$b_{ij} =\big\langle U_i,\ V_j \big\rangle \quad \text{for all} \quad i, j.$$
It remains to notice that  the dimension of the space spanned by vectors $U_i$ and $V_j$ does not exceed 
$$1+\sum_{m=1}^{k} \dim \sym \left({\Bbb R}^d\right)^{\otimes m} =1+\sum_{m=1}^{k} {d+m-1 \choose m} =
{d+k \choose k}.$$
{\hfill \hfill \hfill} \qed
\enddemo

\proclaim{(4.5) Corollary}  Let $A=\left(a_{ij}\right)$ for $i \in I$ and $j \in J$ be a real matrix, let $S \subset {\Bbb R}$,
$$S=\left\{a_{ij}: \quad i \in I,\ j \in J \right\},$$
 be the 
set of all distinct values among the 
matrix entries $a_{ij}$ and let $\phi: S \longrightarrow {\Bbb R}$ be a function.
Let us define a matrix $B=\left(b_{ij} \right)$ by 
$$b_{ij}=\phi\left(a_{ij}\right) \quad \text{for all} \quad i,j.$$
If $|S| \leq k$ then 
$$\rk B \ \leq \ {k-1+ \rk A  \choose k-1}.$$
\endproclaim
\demo{Proof}
Since $|S| \leq k$ there is a
polynomial $p: {\Bbb R} \longrightarrow {\Bbb R}$ with $\deg p \leq k-1$ such that $\phi(t)=p(t)$ for all $t \in S$, so 
$b_{ij}=p\left(a_{ij}\right)$. We write 
$$p(t)=\sum_{m=0}^{k-1} \alpha_m t^m$$
for some $\alpha_m \in {\Bbb R}$. The proof follows by Lemma 4.4.
{\hfill \hfill \hfill} \qed
\enddemo

\proclaim{(4.6) Lemma} Let $A=\left(a_{ij}\right)$ be a real non-negative matrix such that the number of distinct values 
among the matrix entries $a_{ij}$ does not exceed $k$. Then 
$$\rkp A \ \leq \ {k-1+ \rk A  \choose k-1}.$$
\endproclaim
\demo{Proof} Let us define a matrix $B=\left(b_{ij}\right)$ by 
$$b_{ij}=\sqrt{a_{ij}} \quad \text{for all} \quad i, j.$$
By Corollary 4.5, 
$$\rk B \ \leq \ {k-1 +\rk A \choose k-1}.$$
Since 
$$a_{ij}=b_{ij}^2 \quad \text{for all} \quad i, j,$$
the proof follows by Lemma 4.3.
{\hfill \hfill \hfill} \qed
\enddemo

\proclaim{(4.7) Lemma} Let $\left\{u_i: \ i \in I \right\} \subset {\Bbb R}^d$ and $\left\{v_j: \ j \in J \right\} \subset {\Bbb R}^d$ be sets of vectors such that 
$$\langle u_i, v_j \rangle \ \leq \ 1 \quad \text{for all} \quad i, j.$$
Let
$${\Cal L} =\spa\left(v_j: \ j \in J\right).$$
Let us define a matrix $A=\left(a_{ij}\right)$ by 
$$a_{ij} =1 - \langle u_i, v_j \rangle \quad \text{for all} \quad i\in I \quad \text{and all} \quad j \in J.$$
Suppose further, that the number of distinct values among the entries $a_{ij}$ does not exceed some positive integer $k$.

Then there exists a convex set $C \subset {\Bbb R}^d$ such that the following holds:
\roster
\item We can write $C$ as a Minkowski sum $C=C' +{\Cal L}^{\bot}$, where 
$C' \subset {\Cal L}$ is an affine image of a section of the cone of
$r \times r$ symmetric positive semidefinite matrices by an affine subspace
for 
$$r \ \leq \ {k+d \choose k-1};$$
\item We have 
$$u_i \in C \quad \text{for all} \quad i \in I$$
and
\item We have
$$\langle x, v_j \rangle \ \leq \ 1 \quad \text{for all } \quad j \in J  \quad \text{and all} \quad x \in C.$$
\endroster
\endproclaim
\demo{Proof} Since $\rk A \leq d+1$ it follows by Lemma 4.6 that 
$$\rkp A \ \leq \ {k+d \choose k-1}.$$
Let $u_i'$ be the orthogonal projection of $u_i$ onto ${\Cal L}$.
By Lemma 4.2, there exists a convex set $C' \subset {\Cal L}$ 
which is an affine image of a section of the cone of $r \times r$ symmetric positive semidefinite matrices by an affine 
subspace for 
$$r \ \leq \ {k+d \choose k-1},$$
such that 
$$u_i' \in C' \quad \text{for all} \quad i \in I$$
and such that
$$\langle x, v_j \rangle \ \leq 1  \quad \text{for all} \quad x \in C' \quad \text{and all} \quad j \in J.$$
Then $C=C' +{\Cal L}^{\bot}$ satisfies the desired conditions.
{\hfill \hfill \hfill} \qed
\enddemo

\subhead (4.8) Proof of Theorem 1.2 \endsubhead 
We consider ${\Bbb R}^d$ as the coordinate hyperplane $x_{d+1}=0$ of ${\Bbb R}^d$. For a vector $x \in {\Bbb R}^d$ 
and $a \in {\Bbb R}$ we denote by $(x, a) \in {\Bbb R}^{d+1}$ the vector obtained from $x$ by appending the $(d+1)$-st 
coordinate equal to $a$.

For a vector $u \in B$ let $\widehat{u}=(u, 1) \in {\Bbb Z}^{d+1}$ and let 
$\left\{ \widehat{u}_i: \ i \in I \right\}$ be the set of vectors obtained this way.

For any vector $v \in {\Bbb Z}^d$ and any $m \in {\Bbb Z}$ such that 
$$\max_{u \in B} \langle u, v \rangle \ = \ m \quad \text{and} \quad \min_{u \in B} \langle u, v \rangle \ \geq \ m-k $$
we let $\widehat{v}=(v, k-m) \in {\Bbb Z}^{d+1}$ and let 
$\left\{ \widehat{v}_j: \ j \in J \right\}$ be the set of all vectors obtained this way. In particular,
$$\max_{i \in I} \langle \widehat{u}_i,\ \widehat{v}_j \rangle  =k \quad \text{and} \quad \min_{i \in I} \langle \widehat{u}_i,\ 
\widehat{v}_j \rangle \ \geq 0 \quad \text{for all} \quad j \in J.$$
We define matrix $A=\left(a_{ij}\right)$ by 
$$a_{ij} = 1- {1 \over k} \langle \widehat{u}_i, \widehat{v}_j \rangle \quad \text{for all} \quad i \in I \quad \text{and all} \quad j \in J.$$
Hence we have $0 \leq a_{ij} \leq 1$ and $k a_{ij} \in {\Bbb Z}$ for all $i$ and $j$. In particular, the number of 
distinct values among the entries $a_{ij}$ does not exceed $k+1$.  

Let 
$${\Cal L} =\spa\left(\widehat{v}_j: \ j \in J\right).$$
Since vector $(0, k)$ is among vectors $\widehat{v}_j$, we have 
$${\Cal L}^{\bot} \subset {\Bbb R}^d.$$
By Lemma 4.7, there is an affine image $\widehat{C}' \subset {\Cal L}$ of a section of the cone of $r \times r$ symmetric 
positive semidefinite matrices with 
$$ r \ \leq \ {k+d +2 \choose k}$$
such that for $\widehat{C}=\widehat{C}' + {\Cal L}^{\bot}$ we have 
$$\widehat{u}_i \ \in \widehat{C} \quad \text{for all} \quad i \in I \quad \text{and} \quad \langle x, \widehat{v}_j \rangle \ \leq \ k 
\quad \text{for all}
\quad x \in \widehat{C} \quad \text{and all}  \quad j \in J.$$
Let us define 
$$C'=\left\{x \in {\Bbb R}^d: \ (x, 1) \in \widehat{C}' \right\} \quad \text{and} \quad C= C' + {\Cal L}^{\bot}.$$
Since ${\Cal L}^{\bot} \subset {\Bbb R}^d$ we have
$$C=\left\{x \in {\Bbb R}^d:\ (x, 1) \in \widehat{C} \right\}.$$
We observe that $C'$ is an affine image of a section of the cone of $r \times r$ positive semidefinite matrices and that
$$B \subset C.$$
Moreover, if for some $v \in {\Bbb Z}^d$ and some integer $m$ we have 
$$\max_{u \in B}  \langle u, v \rangle \ = \ m  \quad \text{and} \quad \min_{u \in B} \langle u, v \rangle \ \geq \ m-k$$
then $(v, k-m) =\widehat{v}_j$ for some $j \in J$ and
$$\langle x, v \rangle \ \leq \ m \quad \text{for all} \quad x \in C.$$
Since $B \subset C$ we necessarily have 
$$\max_{x \in C} \langle x, v \rangle =m.$$ 
{\hfill \hfill \hfill} \qed

\head 5. Approximating non-negative matrices by matrices of a small positive semidefinite rank \endhead

As a by-product of our proof of Theorem 1.2 we obtain the following result.
\proclaim{(5.1) Theorem} For any $\epsilon > 0$ there is a positive integer $d=d(\epsilon)$ such that the following holds.
Let $A=\left(a_{ij} \right)$, $i \in I$, $j \in J$ be a matrix such that 
$$0 \ \leq \ a_{ij} \ \leq \ 1 \quad \text{for all} \quad i, j.$$
Then there exists a non-negative $m\times n$ matrix $A'=\left(a_{ij}'\right)$ such that 
$$\left| a_{ij} -a_{ij}'\right| \ \leq \ \epsilon \quad \text{for all} \quad i, j$$
and 
$$\rkp A' \ \leq \ \left(  \rk A \right)^d.$$
\endproclaim
\demo{Proof} Without loss of generality, we assume that $\rk A \geq 2$ and that $\epsilon < 1/2$. There is a univariate 
polynomial $p(t)$ of some degree $k=k(\epsilon)$ such that 
$$\left| \sqrt{t} - p(t) \right| \ \leq \ {\epsilon \over 3} \quad \text{for all} \quad 0 \leq t \leq 1. \tag5.1.1$$
We let $B=\left(b_{ij}\right)$ by 
$$b_{ij} =p\left(a_{ij}\right) \quad \text{for all} \quad i, j.$$
By Lemma 4.4 that 
$$\rk B \ \leq \ { k + \rk A \choose k }.$$
Let us define $A'$ by 
$$a_{ij}' = b_{ij}^2 \quad \text{for all} \quad i,j.$$
It follows by (5.1.1) that $a_{ij}'$ approximates $a_{ij}$ within $\epsilon$. 
By Lemma 4.3 
$$\rkp A' \ \leq \ \rk B \ \leq \ {k + \rk A \choose k}$$ 
and the proof follows.
{\hfill \hfill \hfill} \qed
\enddemo

It would be interesting to find out if Theorem 5.1 leads to any non-trivial approximations of general convex bodies 
by projections of spectrahedra. So far, the author was unable to beat the bounds established by Theorem 1.1.

\head Acknowledgment \endhead

The author is grateful to Rekha Thomas for pointing out that papers \cite{G+10} and \cite{GT10} contain the proof of 
Theorem 1.2, though stated in the different language of theta bodies.


\Refs
\widestnumber\key{AAAA}

\ref\key{Ba97}
\by K. Ball
\paper An elementary introduction to modern convex geometry
\inbook Flavors of geometry
\bookinfo Math. Sci. Res. Inst. Publ.
\vol 31
\publ Cambridge Univ. Press
\publaddr Cambridge
\yr1997
\pages 1--58 
\endref

\ref\key{BV08}
\by A. Barvinok and E. Veomett
\paper The computational complexity of convex bodies
\inbook Surveys on discrete and computational geometry
\bookinfo Contemp. Math.
\vol 453
\publ Amer. Math. Soc.
\publaddr Providence, RI
\yr 2008
\pages 117--137
\endref

\ref\key{F+11}
\by S. Fiorini, S. Massar, S. Pokutta, H.R. Tiwary and R. de Wolf
\paper Linear vs. semidefinite extended formulations: exponential separation and strong lower bounds
\paperinfo preprint {\tt arXiv:1111.0837}
\yr 2011
\endref

\ref\key{GN11}
\by J. Gouveia and T. Netzer
\paper Positive polynomials and projections of spectrahedra
\jour SIAM J. Optim. 
\vol 21 
\yr 2011
\pages 960--976
\endref

\ref\key{G+10}
\by J. Gouveia, P.A. Parrilo and R.R. Thomas
\paper Theta bodies for polynomial ideals
\jour SIAM J. Optim.
\vol 20
\pages 2097--2118
\endref

\ref\key{GT10}
\by J. Gouveia and R.R. Thomas
\paper Convex hulls of algebraic sets 
\paperinfo preprint {\tt arXiv:1007.1191}, to appear in ``Handbook of Semidefinite, Cone and Polynomial Optimization:
Theory, Algorithms, Software and Applications"
\yr 2010
\endref

\ref\key{G+11}
\by J. Gouveia, P.A. Parrilo and R.R. Thomas
\paper Lifts of convex sets and cone factorizations
\paperinfo preprint {\tt arXiv:1111.3164}
\yr 2011
\endref

\ref\key{Tu10}
\by L. Tun\c cel
\book Polyhedral and semidefinite programming methods in combinatorial optimization
\bookinfo Fields Institute Monographs
\vol 27
\publ  American Mathematical Society
\publaddr Providence, RI
\yr 2010
\endref

\endRefs
\enddocument
\end